\documentclass[12pt]{amsart}
\usepackage{amssymb}
\newtheorem{thm}{Theorem}[section]
\newtheorem{lem}[thm]{Lemma}
\newtheorem{cor}[thm]{Corollary}
\newtheorem{prop}[thm]{Proposition}
\newtheorem{ex}[thm]{Example}

\newtheorem*{prob*}{Open problem}

\theoremstyle{definition}

\newtheorem{defi}[thm]{Definition}

\theoremstyle{remark}

\newtheorem{rem}[thm]{Remark}
\newtheorem*{rem*}{Remark}

\DeclareMathOperator{\id}{id}

\newcommand{\map}[3]{ #1 : #2 \rightarrow #3 }

\newcommand{\kringel}{\mathbin{\raise1pt\hbox{$\scriptstyle\circ$}}} 
\newcommand{\pkt}{\mathbin{\raise0pt\hbox{$\scriptstyle\bullet$}}}

\newcommand{\C}{\mathbb{C}}

\newcommand{\N}{\mathbb{N}}

\newcommand{\La}{\mathfrak{a}}
\newcommand{\Lb}{\mathfrak{b}}

\newcommand{\Lg}{\mathfrak{g}}
\newcommand{\Lh}{\mathfrak{h}}

\newcommand{\Ll}{\mathfrak{l}}

\newcommand{\Ls}{\mathfrak{s}}

\newcommand{\CM}{\mathcal{M}}

\newcommand{\al}{\alpha}

\newcommand{\Ga}{\Gamma}

\newcommand{\ra}{\rightarrow}

\renewcommand{\phi}{\varphi}

\begin{document}


\title[Faithful representations]{Minimal faithful representations of reductive 
Lie algebras}

\author[D. Burde]{Dietrich Burde}
\author[W. Moens]{Wolfgang Moens}
\address{Fakult\"at f\"ur Mathematik\\
Universit\"at Wien\\
  Nordbergstr. 15\\
  1090 Wien \\
  Austria} 
\email{dietrich.burde@univie.ac.at}
\email {wolfgang.moens@univie.ac.at}
\date{\today}


\begin{abstract}
We prove an explicit formula for the invariant $\mu(\Lg)$ for finite-dimensional 
semisimple, and reductive Lie algebras $\Lg$ over $\C$. 
Here $\mu(\Lg)$ is the minimal dimension of a faithful linear representation of $\Lg$.
The result can be used to study Dynkin's classification of maximal reductive 
subalgebras of semisimple Lie algebras. 
\end{abstract}

\maketitle

\section{Introduction}

In \cite{BU1} the following invariant for finite-dimensional Lie algebras
has been introduced:

\begin{defi}
Let $\Lg$ be an $n$-dimensional Lie algebra over a field $k$. 
Let $\mu(\Lg)$ denote the minimal dimension of a faithful linear
representation of $\Lg$.
\end{defi}

We consider $k$ as given by $\Lg$, so that we need not refer to $k$ in the
notation $\mu(\Lg)$. By Ado's theorem, $\mu(\Lg)$ is finite. 
In general it is not known how to determine this invariant. In particular,
it seems to be very hard in general to determine $\mu(\Lg)$ for a given 
solvable Lie algebra $\Lg$. \\
The invariant $\mu(\Lg)$ plays an important role in the theory of affinely flat
manifolds and affine crystallographic groups, see \cite{BU2}. In particular,
the following two results are known:

\begin{prop}
Let $G$ be an $n$-dimensional Lie group with Lie algebra $\Lg$. If $G$ admits a
left-invaraint affine structure then $\mu(\Lg)\le n+1$.
\end{prop}

\begin{prop}
Let $\Ga$ be a torsionfree finitely generated nilpotent group of rank
$n$ and $G_{\Ga}$ its real Malcev-completion with Lie algebra
$\Lg_{\Ga}$. If $\Ga$ is the fundamental group of a compact complete affine
manifold then $\mu(\Lg_{\Ga})\le n+1$.
\end{prop}

It is well known that a semisimple Lie group $G$ does not admit any left-invariant
affine structures. If $G$ is reductive then the existence problem of left-invariant
affine structures is already quite interesting, see \cite{BU2} and the references
cited there. The existence problem is considered to be hard for solvable and
nilpotent Lie groups. \\
If $\Lg$ has trivial center $Z(\Lg)$, the adjoint representation is faithful and
we have $\mu(\Lg)\le n$. If $\Lg$ is nilpotent, the adjoint representation is not
faithful, and such a result is not even true in general.
Since the classification of representations of nilpotent Lie algebras
is a wild problem, it seems reasonable to expect difficulties in determining
$\mu(\Lg)$. In this case ones tries to obtain good upper and lower
bounds for $\mu(\Lg)$. There is the following result, see \cite{BU2}. 
Let $\Lg$ be a nilpotent Lie algebra of dimension $n$ and nilpotency class
$k$. Denote by $p(j)$ the number of partitions of $j$ and let
$$p(n,k)=\sum_{j=0}^{k}\binom{n-j}{k-j}p(j).$$
Then $\mu(\Lg)\le p(n,k)$.
In particular, with $\al=\frac{113}{40}$, we have
\begin{align*}
\mu(\Lg) & < \frac{\al}{\sqrt{n}}\, 2^n.
\end{align*}
If $\Lg$ is reductive, however, the situation is much better. There are explicit formulas
for $\mu(\Lg)$, in case $\Lg$ is abelian or $\Lg$ is simple.
The aim of this paper is to show the following result:

\begin{thm}\label{1.4}
Let $\Lg$ be a complex reductive Lie algebra and $\Lg=\Ls_1 \oplus \cdots \oplus 
\Ls_{\ell}\oplus \C^{k}$
its decomposition into simple ideals $\Ls_i$ and center $\C^k$. Then the following
formula holds:
\[
\mu(\Lg)=\mu(\Ls_1)+\ldots + \mu (\Ls_{\ell})+ \mu (\C^{k-\ell}),
\]
where $ \mu (\C^{k-\ell})=\lceil \sqrt{2(k-\ell-1)}\rceil $ for 
$k-\ell >1$, $\mu(\C)=1$ and $ \mu (\C^{k-\ell})=0$ for $k-\ell\le 0$.
\end{thm}

\section{Faithful representations}

We start with two simple lemmas.

\begin{lem}\label{2.1}
Let $\Lh$ be a subalgebra of $\Lg$. Then $\mu(\Lh)\le \mu(\Lg)$. Furthermore,
if $\La$ and $\Lb$ are two Lie algebras, then $\mu(\La\oplus \Lb)\le \mu(\La)+\mu(\Lb)$.
\end{lem}

\begin{proof}
The composition of the embedding $\Lh \hookrightarrow \Lg$ and a faithful
representation $\Lg\ra \Lg\Ll_n(\C)$ is a faithful representation of $\Lh$ of degree
$n$. If $\phi$ and $\psi$ are faithful representations of $\La$ respectively $\Lb$,
then $\phi\oplus \psi$ is a faithful representation of $\La\oplus \Lb$.
\end{proof}

\begin{lem}\label{2.2}
Let $\Lg$ be a Lie algebra with trivial center. Then 
$\mu(\Lg \oplus \C)=\mu(\Lg)$. 
\end{lem}

\begin{proof}
We have $\mu(\Lg)\le \mu(\Lg \oplus \C)$. Conversely, let $\rho \colon \Lg\ra
\Lg\Ll_n(\C)$ be a faithful representation of minimal dimension $n=\mu(\Lg)$.
Suppose that there is an $x\in \Lg$ such that $\rho(x)=I_n$ is the identity.
Then, for all $y\in \Lg$,
\begin{align*}
\rho([x,y]) & =[\rho(x),\rho(y)] = [I_n,\rho(y)]=0.
\end{align*}
Since $\rho$ is faithful, we have $x\in Z(\Lg)=0$, which is a contradiction.
It follows that the identity $I_n$ is not in $\rho(\Lg)$, and
hence $\rho (\Lg) \oplus \C\cdot I_n$ yields a faithful representation of 
$\Lg\oplus \C$ of dimension $n$. 
\end{proof}

\subsection{Faithful representations of abelian Lie algebras} 

If $\Lg$ is abelian then there exists an explicit formula for
$\mu(\Lg)$, which only depends on the dimension of $\Lg$. If $V$ is an
$n$-dimensional vector space, then any faithful
representation $\phi \colon \Lg \rightarrow \Lg\Ll(V)$ 
turns $\phi (\Lg)$ into an $n$-dimensional commutative subalgebra 
of the matrix algebra $M_d(K)$. Jacobson \cite{JA2} proved:
\begin{prop}\label{jac}
Let $M$ be a commutative subalgebra of $M_d(K)$ over an arbitrary field
$k$. Then $\dim M\le \lfloor d^2/4 \rfloor +1$ and the bound is attained. 
\end{prop}

For $K=\C$ the result was first proved by I. Schur.
The proposition implies the following result, see \cite{BU1}:
\begin{prop}
Let $\Lg$ be an abelian Lie algebra of dimension $n>1$ over a 
field $k$. Then $\mu(\Lg)= \lceil 2\sqrt{n-1}\rceil $.
\end{prop}

For $n=1$ we have $\mu (\Lg)=1$.

\subsection{Faithful representations of simple Lie algebras}

Let $\Lg$ be a finite-dimensional complex simple Lie algebra. Then every 
non-trivial representation $\rho\colon \Lg \ra \Lg\Ll_n(\C)$
is faithful, since $\ker (\rho)$ is an ideal in $\Lg$. 
Hence any simple $\Lg$-module different from the $1$-dimensional trivial module
$\C$ is faithful. Then $\mu(\Lg)$ is the dimension of the ``smallest'' 
simple module different form the trivial one. Such a minimal
simple module is, with few exceptions, the natural or the adjoint module. 
Their dimensions are well known, see for example \cite{CAR}:
\vspace*{0.5cm}
\begin{center}
\begin{tabular}{c|c|c}
$\Lg$ & $\dim (\Lg)$ & $\mu(\Lg)$ \\
\hline
$A_n, \; n\ge 1$ & $(n+1)^2-1$ & $n+1$ \\
$B_2$ & $10$ & $4$ \\
$B_n, \; n\ge 3$ & $2n^2+n$ & $2n+1$ \\
$C_n, \; n\ge 3$ & $2n^2+n$ & $2n$ \\
$D_n, \; n\ge 4$ & $2n^2-n$ & $2n$ \\
$E_6$ & $78$ & $27$ \\
$E_7$ & $133$ & $56$ \\
$E_8$ & $248$ & $248$ \\
$F_4$ & $52$ & $26$ \\
$G_2$ & $14$ & $7$ \\
\end{tabular}
\end{center}

\subsection{Faithful representations of semisimple Lie algebras}

Assume that $\Lg=\Ls_1\oplus \Ls_2$ with simple Lie algebras $\Ls_1$ and $\Ls_2$.
Let $(\rho_i,V_i)$ be a representation of $\Ls_i$ for $i=1,2$. Then
$\rho_i\circ \pi_i$ is a representation of $\Lg$, where $\pi_i$ is the projection
from $\Lg$ to $\Ls_i$. In the language of $\Lg$-modules this means: 
The $\Ls_1$-module $V_1$ is also a $\Lg$-module via 
$(x,y).v=x.v$ for $x\in \Ls_1, y\in \Ls_2$, and the $\Ls_2$-module $V_2$ is 
also a $\Lg$-module via $(x,y).v=y.v$.
The {\it tensor product} is the $\Lg$-module $V_1\otimes V_2$ defined by,  
for $(x,y)\in \Ls_1\oplus \Ls_2$ and $v_i\in V_i$, 
\begin{align*}
(x,y).(v_1\otimes v_2) & =(x,y).v_1\otimes v_2+v_1\otimes (x,y).v_2 \\
 & = x.v_1\otimes v_2 + v_1\otimes y.v_2
\end{align*}
In the language of representations we denote $V_1\otimes V_2$ also simply by
$\rho_1 \otimes \rho_2$. Writing $z=(x,y)\in \Ls_1\oplus \Ls_2$ this means
\[
(\rho_1 \otimes \rho_2) (z)= \rho_1(x)(v_1)\otimes v_2+ v_1\otimes \rho_2(y)(v_2).
\]

\begin{lem}\label{kernel}
We have $\ker (\rho_1 \otimes \rho_2)= \ker (\rho_1)\oplus \ker (\rho_2)$.
\end{lem}

\begin{proof}
Obviously $\ker (\rho_1)\oplus \ker (\rho_2)\subseteq \ker (\rho_1 \otimes \rho_2)$.
Conversely choose an element $z=(x,y)\in \Ls_1\oplus \Ls_2=\Lg$ which lies in the
kernel of $\rho_1\otimes \rho_2$, i.e., 
\[
\rho_1(x)(v_1)\otimes v_2 + v_1\otimes \rho_2(y)(v_2)=0
\] 
for all $v_i\in V_i$. Using explicit bases for $V_1,V_2$ and $V_1\otimes V_2$
one easily obtains
\[
\rho_1(x)=\al \id_{\mid V_1},\quad \rho_2(x)=-\al \id_{\mid V_2}
\]
with a constant $\al\in \C$. Since $\Ls_i$ is a simple Lie algebra for $i=1,2$ 
and $\rho_1(x), \rho_2(y)$ are traceless linear operators, 
it follows $\al=0$ and $(x,y)\in \ker (\rho_1)\oplus \ker (\rho_2)$.
\end{proof}

We can extend the above easily to the case  $\Lg=\Ls_1\oplus \cdots \oplus \Ls_{\ell}$
with representations  $(\rho_i,V_i)$ for  $i=1,2, \ldots ,\ell$.
We have the following result, see \cite{IWA}:

\begin{thm}
Let $\Lg$ be a semisimple Lie algebra and $\Lg=\Ls_1 \oplus \cdots \oplus \Ls_{\ell}$ 
be a decomposition of $\Lg$ into ideals of $\Lg$. Then every
irreducible representation $(\rho,V)$ of $\Lg$ is equivalent
to the tensor product of $\ell$ irreducible representations 
$(\rho_1\circ \pi_1,V_1), \ldots ,(\rho_{\ell}\circ \pi_{\ell},V_{\ell})$. Conversely, if
$(\rho_i,U_i)$ are arbitrary irreducible representations of $\Ls_i$ for
$i=1,\ldots \ell$, then  
$$(\rho_i\circ\pi_i\otimes \cdots \otimes \rho_{\ell}\circ \pi_{\ell},
U_1\otimes \cdots \otimes U_{\ell})$$ 
is an irreducible representation of $\Lg$.
\end{thm}

Let $\Lg=\Ls_1 \oplus \cdots \oplus \Ls_{\ell}$ be semisimple and $\rho$ be a
representation of $\Lg$. Then, by Weyl's theorem,
\[
\rho=\rho_1\oplus \cdots \oplus \rho_n,
\]
with irreducible representations $\rho_i$ of $\Lg$. 
Each of the $\rho_i$ is the tensor product 
$\rho_i=\rho_{i,1}\otimes \cdots \otimes \rho_{i,\ell}$ 
where $\rho_{i,j}$ is an irreducible representation of $\Ls_i$.
This gives
\begin{align}\label{rho}
\rho=\bigoplus_{i=1}^n \bigotimes_{j=1}^{\ell} \rho_{i,j}
\end{align}

For the dimension of $\rho$ we obtain 
\begin{align}\label{dimrho}
\dim \rho=\sum_{i=1}^n \prod_{j=1}^{\ell} \dim \rho_{i,j}.
\end{align}

\begin{defi}
For a representation $\rho$ of $\Lg$ define the following associated matrix
``of dimensions''
\[
\Phi_{\rho}=(\dim \rho_{i,j})_{i,j}\in M_{n,\ell}(\N).
\]
\end{defi}

\begin{lem}\label{2.8}
Let $\Lg=\Ls_1 \oplus \cdots \oplus \Ls_{\ell}$ be a complex semisimple Lie algebra.
A finite-dimensional representation $\rho$ of $\Lg$ is faithful if and only if the 
matrix $\Phi_{\rho}$ has no column consisting only of $1$'s.
\end{lem}

\begin{proof}
As before, write 
$\rho=\bigoplus_{i=1}^n \bigotimes_{j=1}^{\ell} \rho_{i,j}$. By lemma $\ref{kernel}$
we have
\[
\ker (\rho_i)=\bigoplus_{j=1}^{\ell} \ker (\rho_{i,j})
\]
for $i=1,\ldots n$. Furthermore $\ker (\rho)=\cap_{i=1}^n \ker (\rho_i)$.
If there is a column consisting only of $1$'s, say colum $j$, then $\Ls_j\subset \Lg$
is contained in $\ker (\rho)$, so that $\rho$ is not faithful. Conversely, suppose
that there is no column with only $1$'s. Choose an element $z=\oplus_i z_i \in 
\ker (\rho)$. Fix a coordinate, say $z_j$. Because there is no $1$-column there must
be an $i$ such that $\rho_{i,j}$ is faithful. By assumption we have
$0=\rho_i(z)=\otimes_j \rho_{i,j}(z_j)$. Again by lemma  $\ref{kernel}$ we have
$\rho_{i,j}(z_j)=0$, and hence $z_j=0$. This follows for all $j$, hence $z=0$.
\end{proof}

\begin{prop}\label{2.9}
Let $\Lg=\Ls_1 \oplus \cdots \oplus \Ls_{\ell}$ be a semisimple Lie algebra and
$\Ls_i$ simple ideals of $\Lg$. Then
\[
\mu(\Lg)=\mu(\Ls_1)+\ldots + \mu(\Ls_{\ell}).
\]
\end{prop}

\begin{proof}
Let $\CM$ be the space of all dimension matrices $\Phi_{\rho}$ for faithful representations
$\rho$ of a fixed semisimple Lie algebra $\Lg$. According to \eqref{rho} and \eqref{dimrho}
let $d_{ij}=\dim \rho_{i,j}$. The determination of $\mu(\Lg)$ is equivalent to 
minimizing the function
\[
f\colon \CM \mapsto \N, \quad \Phi_{\rho}\mapsto \sum_{i=1}^n \prod_{j=1}^{\ell} d_{ij}.
\]
By lemma $\ref{2.8}$ no column of a matrix $\Phi_{\rho}\in \CM$ contains only $1$'s.
Denote by $P$ the matrix in $\CM$, which has diagonal elements
$d_{ii}=\mu (\Ls_i)$ and all other elements equal to $1$. Then
\[
f(P)=\sum_{i=1}^{\ell} \mu(\Ls_i).
\]
We will show that this is the minimal value of $f$, i.e., $\mu(\Lg)=f(P)$.
Suppose $D=(d_{ij})\in \CM$ is a 
matrix with minimal value $f(D)$. If there is a row, say row $i$, with more than one
element unequal to $1$, say $d_{ij}$ and $d_{ik}$, then construct a new matrix
$C$, by replacing the $i$th row 
\[
(d_1,\ldots ,d_{ij},d_{j+1},\ldots ,d_{ik},d_{k+1},\ldots ,d_{\ell})
\]
of $D$ by two new rows
\[
\begin{pmatrix}
d_1,\ldots , d_{ij}, d_{j+1},\ldots, 1, d_{k+1},\ldots ,d_{\ell} \\
d_1,\ldots , 1 ,d_{j+1},\ldots, d_{ik} ,d_{k+1},\ldots ,d_{\ell} 
\end{pmatrix}.
\]
Note that the new matrix $C$ really is in $\CM$.
It has one more row than $D$ and satisfies $f(C)\le f(D)$ since 
$a+b\le ab$ for integers $a,b\ge 2$. By assumption $f(C)=f(D)$.
After repeating this finitely many times we arrive at a matrix $B\in \CM$ where 
every row has at most one element different from $1$. 
In fact, a row $(1,\ldots , 1)$ is impossible, because otherwise we remove this
row and obtain still a matrix in $A\in \CM$ with $f(A)<f(D)$, which is a contradiction.
Thus every row of $B$ has a unique entry different from $1$. Similarly it
is impossible that a column of $B$ contains more than one of these unique entries.
This implies that the number of rows and columns of $B$ coincides.
Now $f(B)$ is just the sum of these unique entries.
Because the value $f(B)$ is minimal, the unique entries must correspond to the numbers
$\mu(\Ls_i)$. Hence $f(B)=\sum_{i=1}^{\ell} \mu(\Ls_i)=f(P)$.
\end{proof}

\subsection{Faithful representations of reductive Lie algebras}

Let $\Lg=\Ls_1 \oplus \cdots \oplus \Ls_{\ell}\oplus \C^k$ be a reductive Lie algebra
over $\C$. Denote by $\ell$ the {\it length} of $\Lg$, i.e., the number of simple ideals
$\Ls_i$.

\begin{lem}\label{2.10}
Let $\Lg$ be a reductive Lie algebra of length $\ell$ and center $\C^k$.
Then 
\[
\mu (\Lg)\le \mu ([\Lg,\Lg])+ \mu (\C^{k-\ell}).
\]
\end{lem}

\begin{proof}
If $\ell\le k$ then 
$\Lg=(\Ls_1\oplus \C)\oplus \cdots \oplus (\Ls_{\ell}\oplus \C)\oplus \C^{k-\ell}$
and we obtain 
\begin{align*}
\mu(\Lg) & \le \sum_{i=1}^{\ell} \mu (\Ls_i \oplus \C) + \mu(\C^{k-\ell}) \\
 & =  \sum_{i=1}^{\ell} \mu(\Ls_i)+  \mu(\C^{k-\ell}) \\
 & = \mu (\Ls_1\oplus \cdots \oplus  \Ls_{\ell}) +  \mu(\C^{k-\ell}) \\
 & = \mu ([\Lg,\Lg])+ \mu (\C^{k-\ell}).
\end{align*}
by lemma $\ref{2.1}$, lemma  $\ref{2.2}$ and proposition  $\ref{2.9}$.
If $k\le \ell$ then $\mu(\C^{k-\ell})=0$ and $\Lg$ can be embedded in 
$\Ls_1 \oplus \cdots \oplus \Ls_{\ell}\oplus \C^{\ell}$. Then we have, using
the above argument for $k=\ell$,
\[
\mu(\Lg)\le \mu(\Ls_1 \oplus \cdots \oplus \Ls_{\ell}\oplus \C^{\ell}) \le
\mu([\Lg,\Lg])+ \mu(\C^{k-\ell}). 
\]
\end{proof}

\begin{rem}
The statement of theorem $\ref{1.4}$ is that the inequality of the above lemma is
in fact an equality.
\end{rem}

\begin{defi}
Denote by $C_{\phi}=\{A\in \Lg\Ll_n(\C) \mid [A,\phi(x)]=0 \; \forall \; x\in \Lg \}$
the {\it centralizer} of a Lie algebra representation $\phi\colon \Lg\ra \Lg\Ll_n(\C)$.
\end{defi}

Note that $C_{\phi}$ is a Lie subalgebra of $\Lg\Ll_n(\C)$.

\begin{defi}
A pair of two Lie algebra representations $\phi\colon \Lg_1\ra \Lg\Ll_n(\C)$ and
$\psi\colon \Lg_2\ra \Lg\Ll_n(\C)$ is said to {\it commute}, if
\[
[\phi(x),\psi (y)]=0 \;\forall \; x\in \Lg_1, y\in \Lg_2.
\]
\end{defi}

\begin{lem}\label{split}
Let $\Lg_1,\Lg_2$ be two Lie algebras and suppose that $\Lg_1$ has trivial center.
There is a bijective correspondence between representations as follows:
\begin{itemize}
\item[(1)] A faithful representation $\map{\varphi}{\Lg_1 \oplus \Lg_2}{\Lg\Ll_n(\C)}$ 
induces a pair of commuting representations $(\varphi_1,\varphi_2)$ by inclusion,
given by  
$\map{\varphi_j = \varphi \circ \iota_j}{ \Lg_j }{ \Lg\Ll_n(\C)}$ for $j=1,2$,
where $\iota_j$ is the natural inclusion of $\Lg_j$ into $\Lg_1 \oplus \Lg_2$.
\item[(2)] Conversely a pair of commuting faithful representations 
$\map{\varphi_j}{ \Lg_j }{ \Lg\Ll_n(\C)}$ induces a faithful representation 
$\map{\varphi}{\Lg_1 \oplus \Lg_2}{\Lg\Ll_n(\C)}$ by 
$\varphi = \varphi_1 \circ \pi_1 + \varphi_2 \circ \pi_2$, where
$\pi_j$ is the natural projection of $\Lg_1 \oplus \Lg_2$ onto $\Lg_j$.
\end{itemize}
\end{lem}

\begin{proof}
It is clear that $\phi_1,\phi_2$ are faithful representations. We have
\begin{align*}
[\phi_1(x),\phi_2(y)] & = [\phi(x,0),\phi(0,y)] = \phi([(x,0),(0,y)]) =0.
\end{align*}
This shows $(1)$. For $(2)$, note that $\phi$ is a representation. 
Let $(x,y)\in \ker (\phi)$. This means $\phi_1(x)+\phi_2(y)=0$, so that
\[
\phi_1(x)=-\phi_2(y)\in \phi_1(\Lg_1)\cap \phi_2(\Lg_2)\subseteq Z(\phi_1(\Lg_1))
=Z(\Lg_1)=0.
\]
Since $\ker (\phi_1)=\ker (\phi_2)=0$ we have $(x,y)=(0,0)$, and $\phi$ is faithful.
\end{proof}

Fix a semisimple Lie algebra $\Lg=\Ls_1\oplus \cdots \Ls_{\ell}$ of length
$\ell$, and an integer $n\ge \mu(\Lg)$. We will construct a certain faithful  
representation $\phi \colon \Lg\ra \Lg\Ll_n(\C)$ for each $n\ge \mu(\Lg)$. 
Let $\sigma_i\colon \Ls_i\ra \Lg\Ll_{\mu(\Ls_i)}(\C)$ be faithful representations of 
minimal dimension $\mu(\Ls_i)$ for $i=1,\ldots \ell$.
Denote by $\phi_0$ the one-dimensional trivial representation of $\Lg$, and let
$m\phi_0=\phi_0\oplus \cdots \oplus \phi_0$. 

\begin{defi}
Let $\Lg$ be as above and $m=n-\mu(\Lg)$. Define a representation
$\sigma \colon \Lg\ra \Lg\Ll_n(\C)$ by
\[
\sigma= m\phi_0 \oplus \sigma_1\oplus \cdots \oplus \sigma_{\ell}.
\]
Then $\sigma$ is called the {\it standard block representation of degree $n$} 
for $\Lg$.
\end{defi}

We are interested in determining the centralizer of a faithful representation
of $\Lg$. We have the following result.

\begin{prop}\label{2.15}
Let $\Lg$ be as above and fix an integer $n\ge \mu(\Lg)$. The centralizer of
any faithful representation $\phi\colon \Lg\ra \Lg\Ll_n(\C)$ can be embedded
into the centralizer of the standard block representation of degree $n$ for $\Lg$.
\end{prop}

The proof is split up into three lemmas. Let $\phi\colon \Lg\ra \Lg\Ll_n(\C)$
be a faithful representation. Since centralizers of equivalent representations are 
isomorphic we may assume, by Weyl's theorem, that 
$\phi= \oplus_{j=0}^k m_j\phi_j$ for irreducible, inequivalent representations 
$\phi_j$ of $\Lg$, and some $m_j\in \N$. Again let $\phi_0$ denote the $1$-dimensional
trivial representation. \\
The following lemma is well known, and follows easily from Schur's lemma.
\begin{lem}\label{2.16}
Let $\phi= \oplus_{j=0}^k m_j\phi_j$ as above. Then, as Lie algebras,
\[
C_{\phi}\cong \bigoplus_{j=0}^k \Lg\Ll_{m_j}(\C).
\]
\end{lem}

\begin{cor}
The centralizer of the standard block representation $\sigma$ of degree $n=m+\mu(\Lg)$
is isomorphic to $\Lg\Ll_m (\C)\oplus \C^{\ell}$.
\end{cor}
Denote by $d_j$ the degree of the representation $\phi_j$.
Now associate to $\phi= \oplus_{j=0}^k m_j\phi_j$ the representation
\begin{align}\label{3}
\tilde{\phi}=m\phi_0  \oplus \left( \bigoplus_{j=1}^k \phi_j \right),
\end{align}
so that $\phi$ and $\tilde{\phi}$ have the same degree, equal to $n$.
This means, that $m=m_0+\sum_{j=1}^k (m_j-1)d_j$. Note that $\tilde{\phi}$ is again faithful
by lemma $\ref{2.8}$.

\begin{lem}\label{2.18}
Let $\phi= \oplus_{j=0}^k m_j\phi_j$ as above. Then the centralizer $C_{\phi}$
can be embedded into the centralizer $C_{\tilde{\phi}}$.
\end{lem}

\begin{proof}
By permuting the summands in $\phi$ we may assume that 
$m_1=\ldots = m_r=1$ for some $0\le r \le k$, and $m_j\ge 2$ for $j> r$. 
By lemma $\ref{2.16}$, and using $\Lg\Ll_1(\C)\cong \C$, the centralizer 
$C_{\phi}$ is isomorphic to $\Lg\Ll_{m_0}(\C)\oplus \C^r \oplus \left(
\oplus_{j=r+1}^k \Lg\Ll_{m_j}(\C)\right)$, which can be embedded into
$\Lg\Ll_p(\C)\oplus \C^r$, where $p=m_0+\sum_{j=r+1}^km_j$. On the other hand,
$C_{\tilde{\phi}}\cong \Lg\Ll_m(\C)\oplus \C^k$, where $m=m_0+\sum_{j=r+1}^k (m_j-1)d_j$.
Certainly $\C^r$ can be embedded into $\C^k$ since $r\le k$. To prove the claim
of the lemma it remains to show that $p\le m$. This is true because
of $m_j,d_j\ge 2$ for $j\ge r+1$, so that $d_j\ge \frac{m_j}{m_j-1}$.
\end{proof}

\begin{lem}\label{2.19}
Let $\phi$ and $\tilde{\phi}$ be as above. Then $C_{\tilde{\phi}}$
can be embedded into the centralizer of the standard block representation of degree $n$.
\end{lem}

\begin{proof}
Consider the decomposition \eqref{3} of $\tilde{\phi}$. Then $\rho=\oplus_{j=1}^k\phi_j$ 
is a faithful representation of $\Lg$.
We claim that we can choose $r$ representations $\phi_j$, denoted again by
$\phi_1,\ldots ,\phi_r$, such that their direct sum 
 $\rho'=\oplus_{j=1}^r \phi_j$ is still a faithful representation of $\Lg$, where
$r$ is at most $\ell$, the length of $\Lg$. 
Since $\rho$ is faithful its dimension matrix has no columns consisting only of $1$'s.
Hence for every column $j$ of our $\ell$ columns we may choose a row $i$ 
such that the entry $(i,j)$ is 
different from $1$. To every such row $i$ corresponds a representation $\phi_i$.
Then we have chosen $\ell$ rows, but not necessarily distinct ones.
Pick out the $\phi_i$ for the distinct rows. Their direct sum is a faithful
representation of $\Lg$, since its dimension matrix again has no columns consisting 
only of $1$'s. \\
Now rewrite $\tilde{\phi}$, using $\rho'$, as $\tilde{\phi}=m\phi_0\oplus \rho' \oplus
\left(\oplus_{j=r+1}^k \phi_j \right)$. Comparing dimensions, we have
$n=m+\dim (\rho')+ \sum_{j=r+1}^k d_j$. Note that $r\le k$. Since $\rho'$ is a faithful
representation of $\Lg$ we have $\dim (\rho')\ge \mu(\Lg)$. Using 
$\sum_{j=r+1}^k d_j\ge k-r$ we obtain
\begin{align*}
\mu(\Lg) & \le \dim (\rho') = n-m-\sum_{j=r+1}^k d_j \\
 & \le n-m-k+r.
\end{align*}
Now it follows, also using lemma $\ref{2.16}$ that
\begin{align*}
C_{\tilde{\phi}} & \cong \Lg\Ll_m(\C) \oplus \C^{k-r}\oplus \C^r \\
& \hookrightarrow \Lg\Ll_{m+k-r}(\C)\oplus \C^r \\
& \hookrightarrow  \Lg\Ll_{n-\mu(\Lg)}(\C)\oplus \C^{\ell} \\
& \cong C_{\sigma}. 
\end{align*}
\end{proof}

Now we can prove proposition $\ref{2.15}$: we have $C_{\phi}\hookrightarrow 
C_{\tilde{\phi}} \hookrightarrow C_{\sigma}$ by the two preceding lemmas.

\begin{cor}
Let $\Lg$ be a semisimple Lie algebra as above, and $\La$ be a Lie algebra. 
Then $\Lg\oplus \La$ can be embedded into $\Lg\Ll_n(\C)$ if and only if $\La$ 
can be embedded into $\Lg\Ll_{n-\mu(\Lg)}(\C)\oplus \C^{\ell}$.
\end{cor}

\begin{proof}
Suppose $\La$ can be embedded into $\Lg\Ll_{n-\mu(\Lg)}(\C)\oplus \C^{\ell}\cong C_{\sigma}$.
Then we have a pair of commuting embeddings  $\sigma\colon \Lg \hookrightarrow 
\Lg\Ll_n(\C)$ and $\tau\colon \La\hookrightarrow C_{\sigma}\hookrightarrow \Lg\Ll_n(\C)$.
Lemma $\ref{split}, (2)$ gives an embedding $\Lg\oplus \La\hookrightarrow \Lg\Ll_n(\C)$.
The converse direction follows from part $(1)$ of lemma $\ref{split}$.
\end{proof}

Now we turn to the proof of theorem $\ref{1.4}$. Let 
$\Lg=\Ls_1 \oplus \cdots \oplus \Ls_{\ell}\oplus \C^{k}$ be a complex reductive
Lie algebra. We write $\Lg=\Ls\oplus \La$, where $\Ls$ is semisimple and
$\La=\C^k$. Given any embedding $\Lg\hookrightarrow \Lg\Ll_n(\C)$, the above corollary 
implies that there is an embedding $\La \hookrightarrow \Lg\Ll_{n-\mu(\Lg)}\oplus \C^{\ell}$.
Denote by $\al(\Lg)$ the maximal dimension of a commutative subalgebra of $\Lg$.
We have $\al(\La)=k$ and  
\[
\al(\Lg\Ll_m(\C)\oplus \C^{\ell})=\al(\Lg\Ll_m(\C))+\al(\C^{\ell}) 
= \lfloor m^2/4 \rfloor + 1 +\ell,
\]
since $\al$ is additive, see \cite{MAL}. If $\La$ is a subalgebra of $\Lb$ then 
$\al(\La)\le \al(\Lb)$. It follows that
\[
k\le \lfloor (n-\mu(\Ls))^2/4 \rfloor +1+ \ell.
\] 
This implies $n\ge \lceil 2\sqrt{k-\ell-1}\rceil +\mu(\Ls)= \mu(\C^{k-\ell})+\mu(\Ls)$.
Together with lemma $\ref{2.10}$ the formula of theorem $\ref{1.4}$ follows. \qed

Finally, the following result can be derived from the above corollary in a similar
way.

\begin{prop}
Let $\Ls$ be a semisimple Lie algebra and $\Lg$ be a perfect Lie algebra, i.e., 
satisfying $[\Lg,\Lg]=\Lg$. Then we have $\mu(\Ls\oplus \Lg)=\mu(\Ls)+\mu(\Lg)$.
\end{prop}

\begin{rem}
Theorem $\ref{1.4}$ can be used to classify all reductive subalgebras, 
up to isomorphism, of $\Lg\Ll_n(\C)$.
As an example, for $n=4$ we obtain (note that $\mu(\C^5)=4$)
\begin{align*}
\C^i, i & = 1,\ldots ,5 \\
A_k\oplus \C^{i}, k & = 1,2, 3, \; 0\le i\le 4-k\\
A_1\oplus A_1\oplus \C^i, i & = 0,1,2 \\
C_2\oplus \C^i, i & = 0,1.
\end{align*} 
\end{rem}

\subsection{Maximal reductive subalgebras}

If we have a faithful representation $\phi\colon \Lg \ra \Lg\Ll_n(\C)$ 
then $\phi (\Lg)$ lies in a maximal reductive subalgebra of $\Lg\Ll_n(\C)$.
There is a complete classification of all maximal reductive
Lie subalgebras in semisimple Lie algebras, due to Malcev \cite{MAL}, 
Dynkin \cite{DY1}, \cite{DY2} and Borel \cite{BOS}.
Hence one might wonder if one can use this 
classification to give another proof of theorem $\ref{1.4}$. However it turns
out that this may be quite complicated in general. In some cases however, we can
give a nice, short proof. Consider the following easy example.

\begin{ex}
We have $\mu(A_1\oplus \C^4)=5.$
\end{ex}

In fact, it is obvious that $\Lg = A_1\oplus \C^4$ has a faithful representation
of dimension $5$: the direct sum of the natural representations of 
$A_1\oplus \C=\Lg\Ll_2(\C)$ and $\C^3$. It remains to show that $\Lg$ cannot be 
faithfully embedded into $\Lg\Ll_4(\C)$. Suppose it can, i.e., $\Lg$ is a subalgebra
of $\Lg\Ll_4(\C)$. Denote by $\pi\colon \Lg\Ll_4(\C) \ra \Lg\Ll_4(\C)/Z$ the natural
projection, where $Z$ is the center of $\Lg\Ll_4(\C)$ with $\dim Z=1$. 
Then we claim that $\pi(\Lg)$ is a reductive Lie subalgebra of $A_3$, which is either
isomorphic to $\Lg$, or to $A_1\oplus\C^3$: let $z\in Z\cap \Lg$. Then 
$[x,z]=0$ for all $x\in \Lg$, i.e., $Z\cap \Lg\subseteq Z(\Lg)$ and 
$\dim (Z\cap \Lg)\le 1$. It follows that
$\pi(\Lg)\cong A_1\oplus Z(\Lg)/(Z\cap \Lg)$. If $Z\cap \Lg=0$, then
$\pi(\Lg)\cong \Lg$, otherwise $\dim (Z\cap \Lg)=1$, so that
$\pi(\Lg)\cong A_1\oplus\C^3$.
So $A_1\oplus \C^4$ or $A_1\oplus \C^3$ is a reductive 
subalgebra of $A_3$, hence lies in a maximal one. But these are
exactly the following ones:
\[
C_2, A_2\oplus \C, A_1\oplus A_1,  A_1\oplus A_1\oplus \C.
\]
Here $A_1\oplus A_1$ is not contained in the last one. Denote by $\al(\Lg)$ the
maximal dimension of a commutative subalgebra of $\Lg$. We can compute $\al(\Lg)$
for all reductive Lie algebras. For a new proof of this result of Malcev see the
nice article of Suter \cite{SUT}. We have $\al (A_1\oplus \C^3)=4$ but
$\al(C_2)=\al ( A_2\oplus \C)=\al(A_1\oplus A_1\oplus \C)=3$ and
$\al(A_1\oplus A_1)=2$. If $\Lh_1\subset \Lh_2$ for reductive Lie algebras then 
$\al(\Lh_1)\le \al (\Lh_2)$. It follows that $A_1\oplus \C^3$ and hence also
$\Lg$ cannot be a subalgebra of one of the maximal reductive subalgebras of $A_3$.
This is a contradiction. \qed \\[0.5cm]
In this way one can also prove more generally that 
\[
\mu(A_1\oplus \C^k)=2+\lceil 2\sqrt{k-2} \rceil, \quad k\ge 3.
\]
The following example, however, shows that this method of using Dynkin's results 
will become very complicated in general.

\begin{ex}
Show that $\mu(A_1\oplus C_3\oplus \C^6)=12$.
\end{ex}

Assume that $\Lg=A_1\oplus C_3\oplus \C^6$ could be embedded into 
$\Lg\Ll_{11}(\C)=A_{10}\oplus \C$.
Then we may assume that $A_1\oplus C_3\oplus \C^5$ is a reductive subalgebra of $A_{10}$.
Passing to maximal reductive subalgebras we may assume that  $A_1\oplus C_3\oplus \C^5$
is a reductive subalgebra of one of the following algebras:
\[
A_9\oplus \C, A_1\oplus A_8\oplus \C, A_2\oplus A_7\oplus \C, A_3\oplus A_6\oplus \C, 
A_4\oplus A_5 \oplus \C, B_5.
\]
The invariant $\al$ of these Lie algebras is given by
$26,22,19,17,16,11$ respectively, whereas $\al(A_1\oplus C_3\oplus \C^5)=12$.
Unfortunately, the only possibility which can be excluded immediately then is
$B_5$. Then we have to treat all the other cases, which ramify to even more cases
in the next step, repeating this kind of argument. Moreover, the maximal reductive
subalgebras of other types, different from $A_n$, play a role.


\begin{thebibliography}{99}


\bibitem{BU1} D. Burde: {\it A refinement of Ado's Theorem}. Archiv Math.\
\textbf{70} (1998), 118-127.

\bibitem{BU2} D. Burde: {\it Left-symmetric algebras, or pre-Lie algebras in 
geometry and physics}. Central European J.\ of Math.\ 
\textbf{4}, Nr. 3 (2006), 323-357. 


\bibitem{BOS} A. Borel, J. de Siebenthal, 
{\it Les sous-groupes ferm \'es de rang maximum des groupes de Lie compacts}, 
Comment.\ Math.\ Helv.\ \textbf{23} (1949), 200-221.

\bibitem{CAR} R. Carter: {\it Lie Algebras of Finite and Affine Type}.
Cambridge studies in advanced mathematics \textbf{96} (2005).

\bibitem{DY1} E. B. Dynkin: {\it Semisimple subalgebras of semisimple Lie algebras}.
AMS Transl.\ {\bf 6} (1957),  111-244.

\bibitem{DY2} E. B. Dynkin: {\it Maximal subgroups of classical groups}.
AMS Transl.\ {\bf 6} (1957),  245-379.

\bibitem{IWA} N. Iwahori: {\it On real irreducible representations of Lie algebras}.
Nagoya Math.\ J.\  {\bf 14} (1959), 59-83.

\bibitem{JA2} N. Jacobson: {\it Schur's theorem on commutative matrices}.
Bull.\ Amer.\ Math.\ Soc.\ \textbf{50} (1944), 431-436.

\bibitem{MAL} A. Malcev: {\it On semi-simple subgroups of Lie groups}.
Izvestia Akad.\ Nauk SSSR \textbf{8} (1944), 143-174. 

\bibitem{SUT} R. Suter: {\it Abelian ideals in a Borel subalgebra of a complex simple 
Lie algebra}.  Invent.\ Math.\  \textbf{156}  (2004), 175-221.

\end{thebibliography}
\end{document}